\documentclass{amsart}

\begin{document}

\title{Some binomial series obtained by the WZ-method}

\author{Jes\'{u}s Guillera \\ }

\address{Zaragoza (Spain)}

\email{jguillera@able.es}


\subjclass[2000]{}

\maketitle

\begin{abstract}
Using the WZ-method we find some of the easiest Ramanujan's
formulae and also some new interesting Ramanujan-like sums.
\end{abstract}

\section{The WZ-method}
We recall \cite{zeilberger2} that a discrete function $A(n,k)$ is
hypergeometric or closed form (CF) if the quotients

\begin{equation}\nonumber
{A(n+1,k) \over A(n,k)} \quad {\rm and} \quad {A(n,k+1) \over
A(n,k)}
\end{equation}
are both rational functions.

And a pair of functions $F(n,k)$, $G(n,k)$ is said to be of Wilf
and Zeilberger (WZ) if $F$ and $G$ are closed forms and besides
$$ F(n+1,k)-F(n,k)=G(n,k+1)-G(n,k). $$
In this case H. S. Wilf and D. Zeilberger \cite{zeilberger1} have
proved that there exists a rational function $C(n,k)$ such that
$$ G(n,k)=C(n,k)F(n,k). $$
The rational function $C(n,k)$ is the so-called certificate of the
pair $(F,G)$.

We now define
\begin{equation}\nonumber
H(n,k)=F(n+1,n+k)+G(n,n+k).
\end{equation}
Zeilberger has proved that for every WZ pair $F(n,k)$, $G(n,k)$
the following holds:
\begin{equation}\nonumber
\sum_{n=0}^{\infty} G(n,0)=\sum_{n=0}^{\infty} H(n,0).
\end{equation}

In next sections we use WZ-pairs to get some Ramanujan's formulae
and also some new Ramanujan-like ones.

\section{First WZ-pair}
We consider the following discrete function:
$$G(n,k)={(-1)^n(-1)^k \over 2^{10n}2^{2k}} (20n+2k+3)
{{2k \choose k}^2 {2n \choose n}^2 {4n-2k \choose 2n-k} \over {2n
\choose k}{n+k \choose n}}.$$ The package EKHAD \cite{petkovsek}
allows to obtain the companion
$$F(n,k)=64{(-1)^n(-1)^k \over 2^{10n}2^{2k}}{n^2 \over 4n-2k-1} {{2k \choose k}^2 {2n \choose
n}^2 {4n-2k \choose 2n-k} \over {2n \choose k}{n+k \choose n}}.$$
We derive the result
$$\sum_{n=0}^{\infty} {(-1)^n {4n \choose 2n}{2n \choose n}^2 \over
2^{10n}} (20n+3)={1 \over 2} \sum_{n=0}^{\infty} {{2n \choose n}^3
\over 2^{12n}} (42n+5).$$

We can extend the pair to have sense for every value of $k$, not
only integers, in the following way:
\begin{align}
F(n,k) &= {64 \over \pi^3} {n^2 \over 4n-2k-1} {(-1)^n \cos(\pi k)
\Gamma(2n-k+1/2)\Gamma(n+1/2)^3\Gamma(k+1/2)^2
\over \Gamma(n+k+1)\Gamma(2n+1)^2}, \nonumber \\
G(n,k) &={1 \over \pi^3} (20n+2k+3) {(-1)^n \cos(\pi k)
\Gamma(2n-k+1/2)\Gamma(n+1/2)^3\Gamma(k+1/2)^2 \over
\Gamma(n+k+1)\Gamma(2n+1)^2}. \nonumber
\end{align}
If $k$ is an integer, it is a routine to prove that
$$\sum_{n=0}^{\infty} G(n,k)=\sum_{n=0}^{\infty} G(n,k+1),$$
and this implies applying Carlson's theorem \cite{bailey} that for
every value of $k$, even if $k$ is not an integer,
$\sum_{n=0}^{\infty} G(n,k)=A$, where $A$ is a constant. To
determine the value of the constant, observe that
$$ \lim_{t \to 1/2} \sum_{n=1}^{\infty} G(n,t)=0 \quad \Rightarrow
\quad A = \lim_{t \to 1/2} G(0,t)={8 \over \pi}.$$ And we have
that, independently of the value of $k$,
$$\sum_{n=0}^{\infty}G(n,k)={8 \over \pi}.$$
But then we have also the sum of another family of infinite series
because obviously we immediately get
$$\sum_{n=0}^{\infty}H(n,k)={8 \over \pi}.$$
For $k=0$, we get the following results \cite{ramanujan}:
\begin{align}
\sum_{n=0}^{\infty}G(n,0) &= \sum_{n=0}^{\infty} {(-1)^n {4n
\choose 2n}{2n \choose n}^2 \over 2^{10n}} (20n+3)={8 \over \pi},
\nonumber
\\
\sum_{n=0}^{\infty}H(n,0) &={1 \over 2} \sum_{n=0}^{\infty} {{2n
\choose n}^3 \over 2^{12n}} (42n+5)={8 \over \pi}. \nonumber
\end{align}
For other values of $k$ we obtain also interesting results. For
example, for $k=1/4$ we get
\begin{align}
& {\sqrt 2 \over 8 } \sum_{n=0}^{\infty} {(-1)^n \left( 1 \over 2
\right)_n\left( 1 \over 4 \right)_{2n} \over (n!)^2 \left( 1 \over
4 \right)_ n 2 ^{4n} }  {40n+7 \over 4n+1}  =
{\sqrt {\pi} \over \Gamma({3 \over 4})^2}, \nonumber \\
& {3 \sqrt 2 \over 8 } \sum_{n=0}^{\infty} {\left( 1 \over 2
\right)_{2n}^2 \left( 1 \over 2 \right)_{n} \over (n!)^2 \left( 1
\over 4 \right)_ {2n} \left( {1 \over 4} \right)_n 2^{8n}} {
112n^2+88n+11 \over (8n+1)(8n+5)} = {\sqrt {\pi} \over \Gamma({3
\over 4})^2}. \nonumber
\end{align}

\section{Second WZ-pair}
We consider the following discrete function:
\begin{equation}\nonumber
G(n,k)={(-1)^k \over 2^{16n}2^{4k}} (120n^2+84nk+34n+10k+3) {{2k
\choose k}^3 {2n \choose n}^4 {4n-2k \choose 2n-k} \over {2n
\choose k}{n+k \choose n}^2};
\end{equation}
the package EKHAD \cite{petkovsek} allows to get the companion
\begin{equation}\nonumber
F(n,k)=512{(-1)^k \over 2^{16n}2^{4k}} {n^3 \over 4n-2k-1} {{2k
\choose k}^3 {2n \choose n}^4 {4n-2k \choose 2n-k} \over {2n
\choose k}{n+k \choose n}^2}.
\end{equation}
We have the following result:
\begin{equation}\nonumber
\sum_{n=0}^{\infty} {{4n \choose 2n}{2n \choose n}^4 \over
2^{16n}} (120n^2+34n+3)={1 \over 4} \sum_{n=0}^{\infty} {(-1)^n
{2n \choose n}^5 \over 2^{20n}} (820n^2+180n+13).
\end{equation}
We can extend the pair to have sense for every value of $k$, not
only integers, in the following way:
\begin{multline} \nonumber
F(n,k)={512 \over \pi^5} {n^3 \over 4n-2k-1} \times {\cos(\pi k)
\Gamma(2n-k+1/2)\Gamma(n+1/2)^6\Gamma(k+1/2)^3 \over
\Gamma(n+k+1)^2\Gamma(2n+1)^3},
\end{multline}
\begin{multline} \nonumber
G(n,k)={1 \over \pi^5} (120n^2+84nk+34n+10k+3) \\
\times {\cos(\pi k) \Gamma(2n-k+1/2)\Gamma(n+1/2)^6\Gamma(k+1/2)^3
\over \Gamma(n+k+1)^2\Gamma(2n+1)^3}.
\end{multline}
If $k$ is an integer it is a routine to prove that
\begin{equation}\nonumber
\sum_{n=0}^{\infty} G(n,k)=\sum_{n=0}^{\infty} G(n,k+1),
\end{equation}
and this implies applying Carlson's theorem \cite{bailey} that for
every value of $k$, even if $k$ is not an integer,
$\sum_{n=0}^{\infty} G(n,k)=A$, where $A$ is a constant. To
determine the constant value $A$ observe that
\begin{equation}\nonumber
\lim_{t \to 1/2} \sum_{n=1}^{\infty} G(n,t)=0 \quad \Rightarrow
\quad A=\lim_{t \to 1/2} G(0,t)={32 \over \pi^2}.
\end{equation}
And we have that, independently of the value of $k$,
\begin{equation}\nonumber
\sum_{n=0}^{\infty}G(n,k)={32 \over \pi^2}.
\end{equation}
But then we have also the sum of another family of infinite series
because obviously we immediately get
\begin{equation}\nonumber
\sum_{n=0}^{\infty}H(n,k)={32 \over \pi^2}.
\end{equation}
For $k=0$ we obtain the following results:
\begin{align}
\sum_{n=0}^{\infty}G(n,0) &= \sum_{n=0}^{\infty} {{4n \choose
2n}{2n \choose n}^4 \over 2^{16n}} (120n^2+34n+3)={32 \over
\pi^2}, \nonumber \\
\sum_{n=0}^{\infty}H(n,0) &= {1 \over 4} \sum_{n=0}^{\infty}
{(-1)^n{2n \choose n}^5 \over 2^{20n}} (820n^2+180n+13)={32 \over
\pi^2}. \nonumber
\end{align}

For other values of $k$ we obtain also interesting results. For
example, for $k=1/4$ we get
\begin{align}
& {1 \over 8 } \sum_{n=0}^{\infty} {\left( 1 \over 2 \right)_{n}^3
\left( 1 \over 4 \right)_{2n} \over (n!)^3 \left( 1 \over 4
\right)_ n^2 2 ^{6n}} {240n^2+110n+11 \over (4n+1)^2}=
{\pi \over \Gamma({3 \over 4})^4}, \nonumber \\
& {1 \over 8 } \sum_{n=0}^{\infty} {(-1)^n \left( 1 \over 2
\right)_{2n}^3 \left( 1 \over 2 \right)_{n}^3 \over (n!)^3 \left(
1 \over 4 \right)_ {2n}^2 \left( {1 \over 4} \right)_n^2 2^{12n}}
{26240n^4+41184n^3+21448n^2+4170n+279 \over (8n+1)^2(8n+5)^2} =
{\pi \over \Gamma({3 \over 4})^4}. \nonumber
\end{align}

\enddocument